\documentclass[12pt,a4paper]{article}
\usepackage{mathrsfs}
\usepackage{amssymb}
\usepackage{amsmath}
\usepackage{amsfonts}
\usepackage{amstext}
\usepackage{amsthm}
\usepackage{graphicx}
\usepackage[top=1in, bottom=1in,left=0.8in,right=0.8in]{geometry}

\def \qed {\hfill \vrule height6pt width 6pt depth 0pt}

\begin{document}

\title{Maximum principle and one-sign solutions for the elliptic $p$-Laplacian
\thanks{Research supported by the NSFC (No.11061030).}}
\author{{\small  Guowei Dai\thanks{Corresponding author. Tel: +86 931
7971297.\newline
\text{\quad\,\,\, E-mail address}: daiguowei@nwnu.edu.cn (G. Dai).%
}
} \\
{\small Department of Mathematics, Northwest Normal University, Lanzhou, 730070, PR China}\\
}
\date{}
\maketitle

\begin{abstract}
In this paper, we prove a maximum principle for the $p$-Laplacian with a sign-changing weight. As an application of this maximum principle,
we study the existence of one-sign solutions
for a class of quasilinear elliptic problems.
\end{abstract}\textbf{\ }

\numberwithin{equation}{section}

\numberwithin{equation}{section}

\section{Introduction}

\bigskip

\quad\, Recently, Dai and Ma [\ref{DM}] studied the existence of one-sign solutions
for the following $p$-Laplacian problem
\begin{equation}
\left\{
\begin{array}{l}
-\Delta_pu=\widetilde{f}(x,u)\,\, \text{in\,\,}\Omega,\\
u=0~~~~~~~~~~\quad\,\,\,\text{on}\,\, \partial\Omega,
\end{array}
\right.\nonumber
\end{equation}
where $\Omega$ is a bounded open subset of $\mathbb{R}^N$, $-\Delta_p u=-\text{div}(\vert \nabla u\vert^{p-2}\nabla u)$ with $1<p<N$
is the $p$-Laplacian of $u$, $\widetilde{f}:\Omega \times \mathbb{R}\rightarrow \mathbb{R}$ is a continuous function.
Under the assumption that $\widetilde{f}(x,s)$ satisfies signum condition $\widetilde{f}(x,s)s>0$ for $s\neq 0$ and crosses the first eigenvalue of $-\Delta_p$, they showed that the above problem possesses at least a positive solution and a negative one.

Naturally, one may ask what will happen if $\widetilde{f}$ does not satisfy the signum condition.
The main purpose of this paper is to establish a result similar to that of [\ref{DM}] for the following $p$-Laplacian problem
\begin{equation}
\left\{
\begin{array}{l}
-\Delta_pu=m(x)f(u)\,\, \text{in\,}\,\Omega,\\
u=0~~~~~~~~~~~~~\quad\,\,\,\,\text{on}\,\, \partial\Omega,
\end{array}
\right.
\end{equation}
where $m:\Omega\rightarrow \mathbb{R}$ is a continuous negative function with $m\not\equiv 0$ in $\Omega$ and $f: \mathbb{R}\rightarrow\mathbb{R}$ is continuous and satisfies:
\\

(f1) $f(s)/\varphi_p(s)$ is bounded for $s\in \mathbb{R}\setminus\{0\}$.
\\

It is well-known that
\begin{equation}
\left\{
\begin{array}{l}
-\Delta_pu=\lambda m(x)\varphi_p(u)\,\, \text{in\,}\,\Omega,\\
u=0~~~~~~~~~~~~~~~~\quad\,\,\,\,\text{on}\,\, \partial\Omega
\end{array}
\right.
\end{equation}
possesses one principal eigenvalues $\lambda_1$ (see [\ref{A}]).
Moreover, if $m$ changes sign in $\Omega$, problem (1.2) possesses two principal eigenvalues $\lambda_1^+$ and $\lambda_1^-$ (see [\ref{Cu}]) with
\begin{equation}
\lambda_{1}^+=\inf\left\{\int_\Omega\vert\nabla u\vert^p\,dx\Big| u\in W_0^{1,p}(\Omega), \int_\Omega m\vert u\vert^p\,dx=1\right\}\nonumber
\end{equation}
and
\begin{equation}
\lambda_{1}^-=\max\left\{-\int_\Omega\vert\nabla u\vert^p\,dx\Big| u\in W_0^{1,p}(\Omega), \int_\Omega-m\vert u\vert^p\,dx=1\right\}.\nonumber
\end{equation}

Define
\begin{equation}
f_0=\lim_{s\rightarrow 0}\frac{f(s)}{\varphi_p(s)}\,\,\text{and}\,\,f_\infty=\frac{f(s)}{\varphi_p(s)}.\nonumber
\end{equation}
Furthermore, we also suppose that
\\

(f2) $f_0<\lambda_1<f_\infty$ or $f_0>\lambda_1>f_\infty$.
\\

Obviously, the methods used in [\ref{DM}] cannot be used to deal with problem (1.1) because
we does not require that $f$ satisfies the signum condition $f(s)s>0$ for $s\neq 0$, which raises the essential difficulty. In order to overcome this difficulty, we use a maximum principle which will be proved in Section 2.
More precisely, we consider the following problem
\begin{equation}
\left\{
\begin{array}{l}
-\Delta_pu-\lambda m(x)\varphi_p(u)=h(x)\,\, \text{in\,}\,\Omega,\\
u=0~~~~~~~~~~~~~~~~~~~~~~~~~~\quad\,\,\,\,\text{on}\,\, \partial\Omega,
\end{array}
\right.
\end{equation}
where $m$ is a changing-sign function.\\

Let $p'=p/(p-1)$. The main result of this work is the following maximum principle.
\\ \\
\textbf{Theorem 1.1.} \emph{Let $h\in L^{p'}(\Omega)$ with $h\geq0$ ($\leq 0$), $\not\equiv 0$ in $\Omega$ and $u$ be a solution of problem (1.3). Then $u>0$ ($<0$) in $\Omega$ if and only if $\lambda\in\left(\lambda_1^-,\lambda_1^+\right)$.}
\\

For the case of $p=2$, Hess and Kato [\ref{HK}, Corollary 2] proved that the condition is sufficient.
For $m(x)\equiv 1$ in $\Omega$, Fleckinger et al. [\ref{FHT}, Theorem 5] showed that the condition is sufficient and necessary.
Thus, the result of Theorem 1.1 has extended and improved the corresponding ones to [\ref{FHT}, \ref{HK}].
\\ \\
\indent In particular, we have the following corollary.
\\ \\
\textbf{Corollary 1.1.} \emph{Assume that $m\geq 0$ and $m\not\equiv 0$ in $\Omega$. Let $h\in L^{p'}(\Omega)$ with $h\geq0$ ($\leq 0$), $\not\equiv 0$ in $\Omega$ and $u$ be a solution of problem (1.3). Then $u>0$ ($<0$) in $\Omega$ if and only if $\lambda\in\left(0,\lambda_1\right)$.}
\\ \\
\textbf{Remark 1.1.} Note that the result of Corollary 1.1 is enough for this work. We prove more general result like as Theorem 1.1 because we believe
that it will be useful in dealing with nonlinear problems with indefinite weight. We shall discuss this kind of problems in our future work.
\\

On the basis of Corollary 1.1, we obtain the following result.
\\ \\
\textbf{Theorem 1.2.} \emph{Let (f1) and (f2) hold. Then problem (1.1) possesses at least a positive and one negative solution.}\\
\\
\textbf{Remark 1.2.} Obviously, the result of Theorem 1.2 extends and improves the corresponding one of [\ref{DPM}, Theorem 3.1].
Meanwhile, the result of Theorem 1.2 also extends the results of Theorem 5.1 and 5.2 of [\ref{DM}] in some sense.
\\ \\
\indent From the results of Theorem 1.2, we can easily deduce the following corollary.
\\ \\
\textbf{Corollary 1.2.} \emph{Besides the assumptions of Theorem 1.2, we also suppose that
$f_0$, $f_\infty\in(0,+\infty)$.
Then for each
\begin{equation}
\lambda\in\left(\frac{\lambda_1}{f_0},\frac{\lambda_1}{f_\infty}\right)\cup \left(\frac{\lambda_1}{f_\infty},\frac{\lambda_1}{f_0}\right),\nonumber
\end{equation}
the problem
\begin{equation}
\left\{
\begin{array}{l}
-\Delta_pu=\lambda m(x)f(u)\,\, \text{in\,}\,\Omega,\\
u=0~~~~~~~~~~~~~~\quad\,\,\,\,\text{on}\,\, \partial\Omega
\end{array}
\right.
\end{equation}
possesses at least a positive and one negative solution.}

\section{Proof of Theorem 1.1}

\bigskip

\quad\, In this section, we give the proof of Theorem 1.1. As usually, we use $\Vert\cdot\Vert$ to denotes the
norm of $W_0^{1,p}(\Omega)$.
\\ \\
\textbf{Proof of Theorem 1.1.} We first prove that the condition is necessary. We divide the proof into two cases.

\emph{Case 1}. $\lambda\geq 0$.

If $u$ is a positive solution of problem (1.3) for $h\geq0$, multiplying (1.3) by $u$ and integration over $\Omega$, we obtain
that
\begin{equation}
\int_\Omega \vert \nabla u\vert^p\,dx-\lambda\int_\Omega m\vert u\vert^p\,dx=\int_\Omega h(x)u\,dx\geq0.\nonumber
\end{equation}
It follows that
\begin{equation}
\int_\Omega \vert \nabla u\vert^p\,dx\geq\lambda\int_\Omega m\vert u\vert^p\,dx.\nonumber
\end{equation}
This relationship together with the variational characterization of $\lambda_1^+$ implies that $\lambda\leq \lambda_1^+$.
We claim that $\lambda<\lambda_1^+$. Suppose on the contrary that $\lambda=\lambda_1^+$. Let
$u_1$ be the corresponding eigenfunction to $\lambda_1^+$ with $\left\Vert u_1\right\Vert=1$. Obviously, one has
\begin{equation}
\left\{
\begin{array}{l}
-\Delta_p u_1=\lambda_1^+m \varphi_p\left(u_1\right)\,\, \text{ in}\,\,\Omega,\\
u_1=0~~~~~~~~~~~~~~~~~~~~~\,\,\text{on}\,\,\partial\Omega.
\end{array}
\right.\nonumber
\end{equation}
For any $\varepsilon>0$, we apply Picone's identity [\ref{AH}] to the pair $u_1$, $u+\varepsilon$. We obtain that
\begin{equation}
0\leq\lambda_1^+\int_\Omega mu_1^p\,dx-\int_\Omega\left(\lambda_1^+ m\varphi_p(u+\varepsilon)+h(x)\right)\frac{u_1^p}{(u+\varepsilon)^{p-1}}.\nonumber
\end{equation}
It follows that
\begin{equation}
\lambda_1^+\int_\Omega mu_1^p\,dx<\int_\Omega\left(\lambda_1^+ m\varphi_p(u+\varepsilon)+h(x)\right)\frac{u_1^p}{(u+\varepsilon)^{p-1}}\leq\lambda_1^+\int_\Omega mu_1^p\,dx.\nonumber
\end{equation}
We get a contradiction.

\emph{Case 2}. $\lambda<0$.

We restate problem (1.3) as the following form
\begin{equation}
\left\{
\begin{array}{l}
-\Delta_p u-\widehat{\lambda}\widehat{m}\varphi_p\left(u\right)=h\,\, \text{ in}
\,\,\Omega,\\
u=0~~~~~~~~~~~~~~~~~~~~~~~\,\,\text{on}\,\,\partial\Omega,
\end{array}
\right.\nonumber
\end{equation}
where $\widehat{\lambda}=-\lambda$, $\widehat{m}=-m$.
Let $\left(\widehat{\lambda}_1^+, v_1\right)$ with $v_1>0$ in $\Omega$ and $\left\Vert v_1\right\Vert=1$ be
the corresponding principal eigen-pairs of the problem
\begin{equation}
\left\{
\begin{array}{l}
-\Delta_p u=\widehat{\lambda}\widehat{m}\varphi_p\left(u\right)\,\, \text{ in}
\,\,\Omega,\\
u=0~~~~~~~~~~~~~~~~~~\,\,\text{on}\,\,\partial\Omega.
\end{array}
\right.\nonumber
\end{equation}
By reasoning as above, we obtain $\widehat{\lambda}<\widehat{\lambda}_1^+$. It is well-known that $\widehat{\lambda}_1^+=-\lambda_1^-$.
So we get $\lambda>\lambda_1^-$.

Conversely, assume that $\lambda\in\left(\lambda_1^-,\lambda_1^+\right)$ and $u$ is a solution of problem (1.3) for $h\geq 0$. We also divide the proof into two cases.

\emph{Case 1}. $\lambda\geq 0$.

We obtain by multiplying problem (1.3) by $u^-$ and integration over $\Omega$:
\begin{equation}
\int_\Omega -\Delta_p u u^-\,dx=\lambda\int_\Omega m\varphi_p(u)u^-\,dx+\int_\Omega h(x)u^-\,dx.
\end{equation}
Then, it follows from (2.1) that
\begin{eqnarray}
\lambda_1^+\int_\Omega m\left\vert u^-\right\vert\,dx&\leq&\int_\Omega \left\vert \nabla u^-\right\vert^p\,dx\nonumber\\
&\leq&\lambda\int_\Omega m\left\vert u^-\right\vert^p\,dx-\int_\Omega h(x)u^-\,dx\nonumber\\
&\leq& \lambda\int_\Omega m\left\vert u^-\right\vert^p\,dx.
\end{eqnarray}
So we have
\begin{equation}
\left(\lambda_1^+-\lambda\right)\int_\Omega m\left\vert u^-\right\vert^p\,dx\leq 0\,\,\text{and}\,\,\int_\Omega m\left\vert u^-\right\vert^p\,dx\geq 0.\nonumber
\end{equation}
Hence $u^-\equiv 0$, so that $u\geq 0$.

We rewrite problem (1.3) as the following form
\begin{equation}
\left\{
\begin{array}{l}
-\Delta_p u+\lambda m^-\varphi_p\left(u\right)=h+\lambda m^+\varphi_p(u)\,\, \text{ in}
\,\,\Omega,\\
u=0~~~~~~~~~~~~~~~~~~~~~~~~~~~~~~~~~~~~~~~~~~\,\,\text{on}\,\,\partial\Omega.
\end{array}
\right.\nonumber
\end{equation}
The strong maximum principle of [\ref{Mo}] implies that $u>0$ in $\Omega$.

\emph{Case 2}. $\lambda<0$.

By an argument similar to (2.2), we obtain
\begin{equation}
\left(\lambda_1^--\lambda\right)\int_\Omega m\left\vert u^-\right\vert^p\,dx\leq 0\,\,\text{and}\,\,\int_\Omega m\left\vert u^-\right\vert^p\,dx\leq0.\nonumber
\end{equation}
Thus $u^-\equiv 0$ which follows $u\geq 0$. We rewrite problem (1.3) as
\begin{equation}
\left\{
\begin{array}{l}
-\Delta_p u-\lambda m^+\varphi_p\left(u\right)=f-\lambda m^-\varphi_p(u)\,\ \text{ in}
\,\,\Omega,\\
u=0~~~~~~~~~~~~~~~~~~~~~~~~~~~~~~~~~~~~~~~~~~~\,\,\text{on}\,\,\partial\Omega.
\end{array}
\right.\nonumber
\end{equation}
The strong maximum principle of [\ref{Mo}] implies that $u>0$ in $\Omega$.\qed

\section{Proof of Theorem 1.2}

\quad\, In this section, based on Theorem 1.1 and the bifurcation result of [\ref{DH}], we study the existence of one-sign solutions for
problem (1.1). From now on, for simplicity, we write $X:=W_0^{1,p}(\Omega)$.

Let us define $g(s)=f(s)-f_0\varphi_p(s)$,
then
\begin{equation}
\lim_{s\rightarrow 0}\frac{g(s)}{\varphi_p(s)}=0\,\,\text{and}\,\,\lim_{\vert s\vert\rightarrow+\infty}\frac{g(s)}{\varphi_p(s)}=f_\infty-f_0.\nonumber
\end{equation}
The existence of one- sign solutions of problem (1.1) is related to the following eigenvalue problem
\begin{equation}
\left\{
\begin{array}{l}
-\Delta_pu=\lambda m\varphi_p(u)+m g(u)\,\, \text{in\,}\,\Omega,\\
u=0~~~~~~~~~~~~\,\quad\quad\quad\quad\quad\,\,\,\,\text{on}\,\, \partial\Omega,
\end{array}
\right.
\end{equation}
where $\lambda\in \mathbb{R}$ is a parameter.
Thus showing that problem (1.1) has a solution is equivalent to show that problem (3.1) has a solution for $\lambda=f_0$.

Applying Theorem 4.5 of [\ref{DH}] to problem (3.1), we obtain that
there are two distinct unbounded sub-continua
$\mathscr{C}^{+}$ and $\mathscr{C}^{-}$,
consisting of the continuum $\mathscr{C}$ emanating from $\left(\lambda_1, 0\right)$.
Moreover, $u>0$ in $\Omega$ if $(\lambda,u)\in \mathscr{C}^+\setminus\{0\}$ and $u<0$ if $(\lambda,u)\in \mathscr{C}^{-}\setminus\{0\}$.
\\
\\
 \textbf{Lemma 3.1.} \emph{Let $\sigma\in\{+,-\}$. There is a $C>0$ such that if $(\lambda_*,u)\in \mathscr{C}^\sigma$ then $\vert \lambda_*\vert\leq C$.}
\\ \\
\textbf{Proof.} We have that $\left(\lambda_*,u\right)$ satisfies
\begin{equation}
\left\{
\begin{array}{l}
-\Delta_p u=\left(\lambda_*+\frac{g(u)}{\varphi_p(u)}\right)m\varphi_p(u)\,\, \text{in\,}\,\Omega,\\
u=0~~~~~\quad\quad\quad\quad\quad\quad\quad\quad\quad\,\text{on}\,\, \partial\Omega.
\end{array}
\right.
\end{equation}
We do the proof for the case $u>0$ in $\Omega$, the case $u<0$ being similar.

Problem (3.2) can be written as
\begin{equation}
\left\{
\begin{array}{l}
-\Delta_p u-\left(\lambda_*-\lambda^*\right)m\varphi_p(u)=\left(\lambda^*+\frac{g(u)}{\varphi_p(u)}\right)m\varphi_p(u)\,\, \text{in\,}\,\Omega,\\
u=0~~~~~~~~~~~~~~~~~~~~~~~~~~~~~~~\quad\quad\quad\quad\quad\quad\quad\quad\quad\,\text{on}\,\, \partial\Omega.
\end{array}
\right.\nonumber
\end{equation}
Choosing $\lambda^*$ such that
\begin{equation}
\left(\lambda^*+\frac{g(u)}{\varphi_p(u)}\right)m\geq 0\,\,\text{a.e. in}\,\, \Omega.\nonumber
\end{equation}
It follows from Corollary 1.1 that $0<\lambda_*-\lambda^*<\lambda_1$. Thus, we get
\begin{equation}
\left\vert\lambda_*\right\vert<\max\left\{\left\vert \lambda_1+\lambda^*\right\vert,\left\vert \lambda^*\right\vert\right\}=:C.\nonumber
\end{equation}
This completes the proof.\qed
\\
\\
\textbf{Proof of Theorem 1.2.} We only prove the isolated property of $f_0<\lambda_1<f_\infty$ since the case $f_0>\lambda_1>f_\infty$ is completely analogous. We shall show that $\mathscr{C}^\sigma$ crosses the hyperplane $\left\{f_0\right\}\times X$ in $\mathbb{R}\times X$.
Obviously, we have $f_0<\lambda_1$. Let $\left(\lambda_n, u_n\right) \in \mathscr{C}^\sigma$ where $u_n\not\equiv 0$ satisfies $\left\vert\lambda_n\right\vert+\left\Vert u_n\right\Vert\rightarrow+\infty.$
Lemma 3.1 implies that there exists a positive constant $M$ such that $\left\vert\lambda_n\right\vert\leq M$ for each $n\in \mathbb{N}$.
It follows that $\left\Vert u_n\right\Vert\rightarrow+\infty$ as $n\rightarrow +\infty$.

We divide the equation
\begin{equation}
-\Delta_pu_n-\lambda_n m\varphi_p\left(u_n\right)=mg\left(u_n\right)\nonumber
\end{equation}
by $\left\Vert u_n\right\Vert$ and set $\overline{u}_n = u_n/\left\Vert u_n\right\Vert$.
Theorem 4.5 of [\ref{DH}] implies that $\overline{u}_n>0$ (or $<0$) in $\Omega$.
It follows that $u_n=\overline{u}_n\left\Vert u_n\right\Vert\rightarrow +\infty$ (or $-\infty$) as $n\rightarrow+\infty$. Thus, we have
\begin{equation}
\lim_{n\rightarrow+\infty}\frac{g(u_n)}{\varphi_p(u_n)}=f_\infty-f_0.
\end{equation}
Since $\overline{u}_n$ is bounded in $X$,
after taking a subsequence if
necessary, we have that $\overline{u}_n \rightharpoonup \overline{u}$ for some $\overline{u} \in X$ and $\overline{u}_n
\rightarrow \overline{u}$ in $L^{p'}(\Omega)$.
By (3.3) and the compactness of $R_p:L^{p'}(\Omega)\rightarrow X$ (see [\ref{DPM}, p.229]) we obtain that
\begin{equation}
-\Delta_p\overline{u}=\left(\overline{\lambda}+f_\infty-f_0\right) m\varphi_p\left(\overline{u}\right),\nonumber
\end{equation}
where $\overline{\lambda}=\underset{n\rightarrow+\infty}\lim\lambda_n$, again choosing a subsequence and relabeling it if necessary.

It is clear that $\left\Vert \overline{u}\right\Vert=1$ and $\overline{u}\in \overline{\mathscr{C}^\sigma}\subseteq \mathscr{C}^\sigma$ since $\mathscr{C}^\sigma$
is closed in $\mathbb{R}\times X$.
Therefore $\overline{u}\not\equiv0$, i.e., $\overline{\lambda}+f_\infty-f_0$ is an eigenvalue of problem (1.2).
So we have $\overline{\lambda}=\lambda_1+f_0-f_\infty<f_0$.
Therefore, $\mathscr{C}^\sigma$ crosses the hyperplane $\left\{f_0\right\}\times X$ in $\mathbb{R}\times X$.
\qed

\end{document}